\documentclass[a4paper,10pt]{amsproc}
\usepackage{amsfonts,lingmacros,tree-dvips, amsmath, amssymb, graphicx, mathrsfs}

\usepackage[all]{xy}

\newdir{ >}{!/6pt/@{ }*:(1,0.2)@_{>}*:(1,-0.2)@^{>}}

\theoremstyle{plain}

\newtheorem{theorem}{Theorem}[section]

\newtheorem{proposition}[theorem]{Proposition}
\theoremstyle{definition}
\newtheorem{definition}[theorem]{Definition}

\newtheorem{counter example}[theorem]{Counter Example}

\numberwithin{equation}{section}

\DeclareMathAlphabet{\mathscr}{OT1}{pzc}{m}{it} 
\begin{document}
\Large{
		\title{ON COMPLETELY NON-BAIRE UNION IN CATEGORY BASES}
		
		\author[S. Basu]{Sanjib Basu}
		\address{\large{Department of Mathematics,Bethune College,181 Bidhan Sarani}}
		\email{\large{sanjibbasu08@gmail.com}}
		
		\author[A. Deb Ray]{Atasi Deb Ray}
		\address{\large{Department of Pure Mathematics, University of Calcutta, 35, Ballygunge Circular Road, Kolkata 700019, West Bengal, India}}
		\email{\large{debrayatasi@gmail.com}}
		
		\author[A.C.Pramanik]{Abhit Chandra Pramanik}
		\address{\large{Department of Pure Mathematics, University of Calcutta, 35, Ballygunge Circular Road, Kolkata 700019, West Bengal, India}}
		\email{\large{abhit.pramanik@gmail.com}} 
		
		\thanks{The third author thanks the CSIR, New Delhi – 110001, India, for financial support}
	\begin{abstract}
In this paper, we intend to show that under not too restrictive conditions, results much stronger than the one obtained earlier by Hejduk could be established in category bases.
	\end{abstract}
\thanks{}
	\maketitle

\section{INTRODUCTION}
Hejduk [2] proved the following theorem 
\begin{theorem}
	Let (X,$\mathcal{S}$) be any arbitrary category base and $\mathcal{M}$($\mathcal{S}$) be the $\sigma$-ideal of all meager sets in the base (X,$\mathcal{S}$), satisfying the following conditions :
	\begin{enumerate}
		\item For an arbitrary cardinal $\alpha < car(X)$, the family $\mathcal{M}$($\mathcal{S}$) is $\alpha$-additive, i.e., this family is closed under the union of arbitrary $\alpha$-sequence of sets.
		\item There exists a base P of $\mathcal{M}$($\mathcal{S}$) of cardinality not greater than that of X.
	\end{enumerate}
Thus if X$\notin$ $\mathcal{M}$($\mathcal{S}$), then for an arbitrary family \{X$_t$\}$_{t\in T}$ of meager sets, being a partition of X, there exists a set T$'$$\subseteq$T such that $\bigcup\limits_{t\in T'}{X_t}$ is not a Baire set.
\end{theorem}

For our purpose, we need the following Definitions and Theorems. For references, the reader may see [3].
\begin{definition}
A category base is a pair (X,$\mathcal{C}$) where X is a non-empty set and $\mathcal{C}$ is a family of subsets of X, called regions satisfying the following set of axioms:
\begin{enumerate}
\item Every point of X belongs to some region; i,e., X=$\cup$$\mathcal{C}$.
\item Let A be a region and $\mathcal{D}$ be a non-empty family of disjoint regions having cardinality less than the cardinality of $\mathcal{C}$.\\
i) If A$\cap$($\cup$$\mathcal{D}$) contains a region, then there is a region D$\in$$\mathcal{D}$ such that A$\cap$D contains a region.\\
ii) If A$\cap$($\cup$$\mathcal{D}$) contains no region, then there is a region B$\subseteq$A that is disjoint from every region in $\mathcal{D}$.
\end{enumerate}
\end{definition}
\begin{definition}
In a category base (X,$\mathcal{C}$), a set is called singular if every region contains a subregion which is disjoint from the set. Any set which can be expressed as countable union of singular sets is called meager. We denote by $\mathcal{M}$ the class of all meager sets in (X,$\mathcal{C}$) which is obviously a $\sigma$-ideal. A set which is not meager is called abundant. A category base (X,$\mathcal{C}$) is called Baire base if every member of $\mathcal{C}$ is an abundant set.
\end{definition}
\begin{definition}
	In a category base (X,$\mathcal{C}$), a set S is called Baire if in every region, there is a subregion in which either S or its complement X$-$S is meager. We denote by $\mathcal{B}$ the class of all Baire sets in (X,$\mathcal{C}$) which is obviously a $\sigma$-algebra.
\end{definition}
\begin{theorem}
	The intersection of two regions either contains a region or is a singular set. 
\end{theorem}
\begin{proposition}
	 If (X, $\mathcal{C}$) is a category base, $\mathcal{N}$ is a subfamily of $\mathcal{C}$ with the property that each region in $\mathcal{C}$ contains a region in $\mathcal{N}$ and Y= $\cup$ $\mathcal{N}$, then (Y, $\mathcal{N}$) is also a category base and the $\mathcal{N}$-singular sets coincide with the $\mathcal{C}$-singular subsets of Y. In addition, if U is a subset of Y and Y $-$ U is $\mathcal{N}$-singular, then X $-$ U is $\mathcal{C}$-singular.
\end{proposition}
\begin{proposition}
	If (Y,$\mathcal{N}$) is a category base, then there exists a disjoint family $\mathcal{M}$ of $\mathcal{N}$ such that Y $-$ $\cup$$\mathcal{M}$ is a singular set. Moreover, $\mathcal{M}$ can be so selected that for every region N, there exists a region M $\in$ $\mathcal{M}$ such that N $\cap$ M contains a region.
\end{proposition}
\begin{theorem}(The Fundamental Theorem)
	Every abundant set in a category base (X, $\mathcal{C}$) is abundant everywhere in some region. This means that for any abundant set A, there exists a region C in every subregion D of which A is abundant.\\
\end{theorem}

To the above Definitions and Theorems, we further add that in any category base (X,$\mathcal{C}$),
\begin{definition}
	A subfamily $\mathcal{C'}$ of $\mathcal{C}$ is called a $\pi$-base [2] if every region D $\in$ $\mathcal{C}$ contains a region E $\in$ $\mathcal{C'}$.
\end{definition}

\begin{definition}
A set A is called completely non-Baire in a region D if for every B $\in$ $\mathcal{B}$ such that B $\cap$ D is abundant, both A $\cap$ B and (D $-$ A) $\cap$ B are abundant. A set which is completely non-Baire in every region is called completely non-Baire.
\end{definition}

The above definition is analogous to the notion  of a `completely I-nonmeasurable set' given in [4].
\section{RESULTS}
To start with, we make the following assumptions :\\
\begin{enumerate}

\item[(i)] The cardinality \#(X) of X is regular.\\
	
\item[(ii)] Our category base (X,$\mathcal{C}$) is a Baire base and possesses a $\pi$-base $\mathcal{C'}$ with \#($\mathcal{C'}$) not exceeding \#(X).\\
	
\item[(iii)] The system (X,$\mathcal{B}$,$\mathcal{M}$) is \#(X)-additive which means that  $\cup$ $\mathcal{E}$ $\in$ $\mathcal{M}$ whenever $\mathcal{E}$ $\subseteq$ $\mathcal{M}$ and \#($\mathcal{E}$) $<$ \#(X) [1].\\
	
\item[(iv)] \#(X) = Min \{$\chi$ : $\mathcal{C}$ is $\chi$-saturated\} where by the phrase ``$\mathcal{C}$ is $\chi$-saturated'' we mean that if $\mathcal{E}$ $\subseteq$ $\mathcal{C}$ such that \#($\mathcal{E}$) = $\chi$, there are two distinct members E, F $\in$ $\mathcal{E}$ such that E $\cap$ F $\neq$ $\phi$.\\
\end{enumerate}

We also use the combinatorial theorem [2] stated below
\begin{theorem}
	If X is an infinte set and $\Phi_1$ is a family of subsets of X such that 
	\begin{enumerate}
		\item[(i)] \#($\Phi_1$) $\leq$ \#(X).
		\item[(ii)] \#(Z) = \#(X) for all Z $\in \Phi_1$.
	\end{enumerate}
then there exists a family $\Phi_2$ of subsets of X such that
\begin{enumerate}
	\item[(a)] \#($\Phi_2$) $>$ \#(X).
	\item[(b)] Z$_1 \neq Z_2$ implies that \#(Z$_1 \cap Z_2$) $<$ \#(X) for all Z$_1, Z_2 \in \Phi_2$.
	\item[(c)] \#(Z $\cap$ Y) = \#(X) for all Y $\in$ $\Phi$$_1$, Z $\in \Phi_2$.
\end{enumerate}
\end{theorem}
Let $\mathcal{K}$ denote the family of all sets whose complements are members of $\mathcal{C'}$. Let $\bigcap \limits_{< \#(X)}\mathcal{K}$ : all sets representable as intersection of subfamilies of $\mathcal{K}$ whose cardinality is less than \#(X).\\

$\bigcup\limits_{\sigma}\bigcap\limits_{< \#(X)}\mathcal{K}$ : all sets  representable as countable union of sets from $\bigcap \limits_{< \#(X)}\mathcal{K}$.
\begin{theorem}
	Any singular set is contained in a singular set which belongs to the family $\bigcap \limits_{< \#(X)}\mathcal{K}$. Any set in $\mathcal{M}$ is a subset of some set in $\mathcal{M}$ $\cap$ ( $\bigcup\limits_{\sigma}\bigcap\limits_{< \#(X)}\mathcal{K}$).
	\begin{proof}
		Let A be any singular set in (X, $\mathcal{C}$). Since $\mathcal{C'}$ is a $\pi$-base, for every C $\in$ $\mathcal{C}$, there is a set D $\in$ $\mathcal{C'}$ such that D is disjoint from A. This constitutes a subfamily $\mathcal{N}$ of $\mathcal{C}$ satisfying the above property such that every member of $\mathcal{C}$ contains a member of $\mathcal{N}$. Let Y = $\cup$$\mathcal{N}$. Then (Y, $\mathcal{N}$) is a category base and by Proposition 1.7, a subfamily $\mathcal{M}$ of $\mathcal{N}$ can be so selected that for every N $\in$ $\mathcal{N}$, there exists M $\in$ $\mathcal{M}$ such that N $\cap$ M contains a region. This makes Y$-$($\cup$ $\mathcal{M})$ $\mathcal{N}$-singular and consequently by Proposition 1.6, X$-$($\cup$$\mathcal{M}$) is $\mathcal{C}$-singular or singular. Moreover, A $\subseteq$ X$-$($\cup$$\mathcal{M}$).\\
		
		Now by our assumption (iv), X$-$($\cup$$\mathcal{M}$) belongs to the family $\bigcap \limits_{< \#(X)}\mathcal{K}$. Hence A is contained in a singular set which belongs to the family $\bigcap \limits_{< \#(X)}\mathcal{K}$. Therefore, by definition, any member of $\mathcal{M}$ is a subset of some set in $\mathcal{M}$ $\cap$ ( $\bigcup\limits_{\sigma}\bigcap\limits_{< \#(X)}\mathcal{K}$).
		\end{proof}
\end{theorem}
Since \#(X) is a regular cardinal (by assumption (i)), \#($\mathcal{M}$ $\cap$ (  $\bigcup\limits_{\sigma}\bigcap\limits_{< \#(X)}\mathcal{K}$))) does not exceed \#(X).
\begin{theorem}
	Every abundant Baire set, i.e., every set in the family  $\mathcal{B} - \mathcal{M}$ contains a set of the form F $-$ G where F $\in$ $\mathcal{C'}$ and G $\in$ $\mathcal{M}$ $\cap$ (  $\bigcup\limits_{\sigma}\bigcap\limits_{< \#(X)}\mathcal{K}$).
	\begin{proof}
		Let S $\in$ $\mathcal{B} - \mathcal{M}$. Then there exists C $\in$ $\mathcal{C}$ such that C $-$ S $\in$ $\mathcal{M}$. By theorem 2.1, C $-$ S $\subseteq$ G $\in$ $\mathcal{M}$ $\cap$ (  $\bigcup\limits_{\sigma}\bigcap\limits_{< \#(X)}\mathcal{K}$) . Let F $\in$ $\mathcal{C'}$ such that F $\subseteq$ C. Choose the set F $-$ G which proves the theorem.
		\end{proof}
\end{theorem}

Let $\mathcal{H}$ =
\{F $-$ G : F $\in$ $\mathcal{C'}$, G $\in$  $\mathcal{M}$ $\cap$ (  $\bigcup\limits_{\sigma}\bigcap\limits_{< \#(X)}\mathcal{K}$) \}.
From the definition of the class, it follows that \#($\mathcal{H}$) $\leq$ \#(X) and for any H $\in$ $\mathcal{H}$, let Y(H) = \{t $\in$ T : X$_t\cap H \neq \phi$\}.
Since the system (X, $\mathcal{B}$, $\mathcal{M}$) is \#(X)-additive (by assumption (iii)) and the family \{X$_t\}_{t \in T}$ gives rise to a partition of X, so \#(Y(H)) = \#(X).\\

Let $\Phi_1$ = \{Y(H) : H $\in$ $\mathcal{H}$\}. Then \#($\Phi_1$) $\leq$ \#(X). We now use the combinatorial theorem (Theorem 2.1) stated above by virtue of which there exists a family $\Phi_2$ of subsets of T satisfying the properties stated in Theorem 2.1.\\

For any Z $\in$ $\Phi_2$, let X(Z) = $\cup$ \{X$_t : t \in Z\}$. We claim that X(Z) is completely non-Baire. Since every set from the family $\mathcal{B}$ $-$ $\mathcal{M}$ contains a set which belongs to the family $\mathcal{H}$, the claim will be settled if we can show that for any Z $\in$ $\Phi_2$ and any H $\in$ $\mathcal{H}$, both H $\cap$ X(Z) and H $-$ X(Z) are abundant. The first one is obvious since \#(Z $\cap$ Y) = \#(X) for all Y $\in$ $\Phi_1$ and Z $\in$ $\Phi_2$ ( (c) of Theorem 2.1) and the second one follows from (b) of Theorem 2.1 and the fact that (X, $\mathcal{B}$, $\mathcal{M}$) is \#(X)-additive (assumption (iii)). Thus 
\begin{theorem}
	If conditions (i), (ii), (iii) and (iv) are satisfied for any category base (X, $\mathcal{C}$) and X = $\bigcup\limits_{t\in T}$$ X_t$ is a partition of X into meager sets, then there exists a family $\Phi_2$ of subsets of T such that \#($\Phi_2$) $>$ \#(X) and for every Z $\in$ $\Phi_2$, X(Z) is completely non-Baire. Moreover, for $Z_1, Z_2 \in \Phi_2 (Z_1 \neq Z_2), X(Z_1)
	\cap X(Z_2) \in \mathcal{M}$.
\end{theorem}	
	The proof of the last part of the above theorem follows from (b) of Theorem 2.1 and assumption (iii).\\

If we prefer to call two sets A and B as almost disjoint in a category base if A $\cap$ B is meager, then what the above theorem indicates is this that under certain circumstances, it is possible to obtain from a partition of X by meager sets more than \#(X) number of completely non-Baire unions which are mutually almost disjoint.\\

In the next theorem, we show that even if \{X$_t\}_{t\in T}$ is a family of disjoint sets but not a partition of X, still it is possible to obtain at least \#(X) number of mutually disjoint completely non-Baire unions by slightly changing the hypothesis. Here we have adapted the technique used to prove Theorem 2.1 [4].
\begin{theorem}
	If conditions (i), (ii), (iii) and (iv) are satisfied for any category base (X, $\mathcal{C}$) and \{X$_t\}_{t\in T}$ is a family of disjoint meager sets such that for every C $\in$ $\mathcal{C}$, ($\bigcup\limits_{t\in T}$X$_t$) $\cap$ C is abundant, then there are at least \#(X) number of completely non-Baire unions which are mutually disjoint.
	\begin{proof}
		Without loss of generality, we may assume that \#($\mathcal{H}$) = \#(X). Also as ($\bigcup\limits_{t\in T}$X$_t$) $\cap$ C is abundant for every C $\in$ $\mathcal{C}$, so for any H $\in$ $\mathcal{H}$, \#(\{t $\in$ T : X$_t \cap H \neq \phi\}$) = \#(X) by condition (iii). Let $\alpha$ denote the smallest ordinal representing \#(X) and we enumerate all the set in $\mathcal{H}$ = \{H$_\beta : \beta < \alpha \}$. By transfinite induction we will construct a sequence \{(X$_{\xi,\eta}, d_\xi) \in \{X_t\}_{t\in T} \times H_\xi : \xi, \eta < \alpha$\} such that 
		\begin{enumerate}
			\item X$_{\xi,\eta} \cap H_{\xi} \neq \phi$ for all $\xi$,$\eta$ $<$ $\alpha$.
			\item $\bigcup\limits_{\xi,\eta < \alpha}X_{\xi,\eta} \cap \{d_\xi : \xi < \alpha$\} = $\phi$.
			\item X$_{\xi,\eta}$ $\neq$ $X_{\xi^{\prime}, \eta^{\prime}}$ for all $\xi$, $\xi^{\prime}$ $<$ $\alpha$ and for all $\eta$, $\eta^{\prime}$ $<$ $\alpha$ ($\eta$ $\neq$ $\eta^{\prime}$).
				\end{enumerate}
			
			We fix $\beta$ $<$ $\alpha$ and suppose that we have already defined the sequence \{(X$_{\xi,\eta}, d_\xi) \in \{X_t\}_{t\in T} \times H_\xi : \xi, \eta < \beta$\} satisfying the conditions :
			\begin{enumerate}
				\item  X$_{\xi,\eta} \cap H_{\xi} \neq \phi$ for all $\xi$, $\eta$  $<$ $\beta$.
				\item   $\bigcup\limits_{\xi,\eta < \beta}X_{\xi,\eta} \cap \{d_\xi : \xi < \beta$\} = $\phi$.
				\item X$_{\xi,\eta}$ $\neq$ $X_{\xi^{\prime}, \eta^{\prime}}$ for all $\xi$, $\xi^{\prime}$ $<$ $\beta$ and for all $\eta$, $\eta^{\prime}$ $<$ $\beta$ ($\eta$ $\neq$ $\eta^{\prime}$).
			\end{enumerate}
		Let X($d_\xi$) denote that unique X$_t$ such that $d_\xi$ $\in$ X$_t$ and choose X$_{\beta,\eta} \in \{X_t\}_{t\in T}$ such that 
		\begin{enumerate}
			\item X$_{\beta,\eta} \neq X_{\beta,\eta^{\prime}}$ for all $\eta$, $\eta^{\prime} < \beta (\eta \neq \eta^{\prime})$.
			\item X$_{\beta,\eta} \cap H_\beta \neq \phi$ for all $\eta < \beta$.
			\item $\bigcup\limits_{\eta < \beta} X_{\beta,\eta} \cap (\cup \{X(d_\xi) : \xi < \beta\}) = \phi$.
		\end{enumerate}
	
	The choice is justified by virtue of the fact stated at the beginning of the proof. Also, owing to the same reasoning, we can further choose X$_{\xi,\beta} \in \{X_t\} (\xi \leq \beta)$ and $d_\beta \in H_\beta$ such that 
	\begin{enumerate}
		\item X$_{\xi,\beta} \neq X_{\xi^{\prime}, \beta}$ for all $\xi$, $\xi^{\prime} \leq \beta$.
		\item X$_{\xi,\beta} \cap H_\xi \neq \phi$ for all $\xi$ $\leq$ $\beta$.
		\item X$_{\xi,\beta} \cap (\cup \{X(d_{\xi ^{\prime}}) : \xi^{\prime} < \beta\}) = \phi$ for all $\xi$ $\leq$ $\beta$.
		\item ($\bigcup\limits_{\xi,\eta \leq \beta}X_{\xi,\eta}) \cap \{d_\beta\} = \phi$.
	\end{enumerate}

It is not hard to check that each member in the family \{$\bigcup\limits_{\xi < \alpha}X_{\xi,\eta} : \eta < \alpha\}$ of mutually disjoint sets is a completely non-Baire union.
		\end{proof}
\end{theorem}

In any Baire base (X, $\mathcal{C}$), ($\cup X_t) \cap C$ is abundant for every C $\in$ $\mathcal{C}$ provided X = $\bigcup\limits_{t\in T}X_t$ is a partition of X. So by the above theorem in any category base (X, $\mathcal{C}$) if conditions (i)$-$(iv) are satisfied, then given a partition X = $\bigcup\limits_{t\in T}X_t$ of X into a family of meager sets, it is possible to construct at least \#(X) number of completely non-Baire unions which are mutually disjoint. In Theorem 2.4, we have theoretically established the existence of more than \#(X) number of completely non-Baire unions which are mutually almost disjoint. But we are not certain whether they are mutually disjoint. However, we can construct at least \#(X) number of completely non-Baire unions which are mutually disjoint.

\bibliographystyle{plain}

	\end{document}